\documentclass[11pt]{article}
\usepackage[margin=2.5cm]{geometry}
\usepackage{amsmath, amsthm, amsfonts, amssymb}
\usepackage[mathscr]{euscript}
\usepackage{amsmath}
\usepackage{array, float}

\usepackage{graphicx, latexsym}
\usepackage{lscape}
\usepackage{multicol}
\usepackage{multirow}
\usepackage{graphicx}
\usepackage{subfig}
\usepackage{tikz}
\usetikzlibrary{calc,positioning}
\usepackage{cite}
\newtheorem{thm}{Theorem}

\newtheorem{cor}{Corollary}

\newtheorem{lem}[thm]{Lemma}
\newtheorem{prop}[thm]{Proposition}

\newtheorem*{rem}{Remark}

\newtheorem{discu}{Discussion:}

\newtheorem{conje}{Conjecture:}

\usepackage{fancyhdr}
\tikzset{
    vertex/.style={
        circle,
        draw=black,
        fill=white,
        line width=0.55pt,
        minimum size=3mm,
        inner sep=0pt
    },
    redvertex/.style={
        vertex,
        fill=red
    },
    bluevertex/.style={
        vertex,
        fill=blue!75!black
    },
    edge/.style={
        draw=black,
        line width=0.75pt
    },
    rededge/.style={
        draw=red!75!black,
        line width=0.75pt
    },
    blueedge/.style={
        draw=blue!65!black,
        line width=0.75pt
    },
    reddash/.style={
        draw=red!75!black,
        line width=0.75pt,
        dash pattern=on 5pt off 3pt
    },
    bluedash/.style={
        draw=blue!65!black,
        line width=0.75pt,
        dash pattern=on 5pt off 3pt
    },
    blackdash/.style={
        draw=black,
        line width=0.75pt,
        dash pattern=on 5pt off 3pt
    }
}

\begin{document}
	\begin{center}
		\textbf{\huge Geodesic Transversal of Hypercube and Butterfly networks}\\
	\end{center}
	\begin{center}
		\textbf{R. Prabha$^{\text a}$, B. Sandhiya$^{\text a,*}$, Paul Manuel$^{\text b}$}
	\end{center}
	\begin{center}
		\small $^{\text a}$Department of Mathematics, Ethiraj College for Women, Chennai 600 008, India\\
		\small $^{\text b}$Department of Information Science, College of Computing Science and Engineering, Kuwait University, Kuwait\\
	\end{center}
	\begin{center}
		*Corresponding author: sandhiyabasskar14@gmail.com
	\end{center}
\begin{abstract}
    A geodesic $P$ of a connected graph $G(V,E)$ is called a maximal geodesic if there exists no geodesic $P'$ in $G$ that contains $P$ as its sub-path. The geodesic transversal problem determines a subset $S\subseteq V(G)$ such that every possible maximal geodesic of $G$ includes at least one vertex of $S$. Such a set with minimum number of elements is called a gt-set and its cardinality is the gt-number of $G$. In this paper, we  establish a general result, called the \textit{Geodesic Transversal lemma} which provides a lower bound for the gt-number of an arbitrary graph $G$ in terms of the gt-number of its induced subgraphs. Using this lemma, we determine the gt-number of the Butterfly and Hypercube networks. 
\end{abstract}
\hspace*{0.5cm}\textbf

 \section{Introduction}
 Let $G$ be a connected graph with vertex set $V(G)$ and edge set $E(G)$. For $x,y\in V(G)$, the shortest path between $x$ and $y$ is called a \textit{geodesic}. A geodesic is said to be a \textit{maximal geodesic} if it is not a sub-path of any other geodesic of $G$. In a subgraph $H$ of $G$, a maximal geodesic $P(u,v)$ in $H$ either remains as a maximal geodesic in $G$ or is contained in another geodesic $P'(x,y)$ in $G$. This geodesic $P'$ is called an extended geodesic and if it is maximal, we call it an extended maximal geodesic of $G$. A \textit{geodesic transversal}  is a set $S$ of vertices of $G$ such that every maximal geodesic in $G$ contains at least one vertex from $S$. A geodesic transversal with minimum cardinality is called the \textit{geodesic transversal set} or \textit{gt-set} of $G$ and its cardinality is called the \textit{geodesic transversal number} or \textit{gt-number} of $G$, denoted as $gt(G)$.

 The geodesic transversal problem was initially proposed by Bre\v{s}ar, Klav\v{z}ar and Manuel \cite{geo.trans} in 2022. The paper illustrates the connection between the geodesic transversal problem and a broad class of existing hitting set problems. The complexity of the gt-problem is proved to be NP-complete, though for trees and spread cacti, it can be computed in polynomial time \cite{geo.trans}. In \cite{geo.trans.net}, the authors determine the gt-number of certain networks like mesh of trees, silicate network and $r$-layer triangular carbon nanosheet using the concept of geodesic packing of graphs~\cite{geo.pack}.

 The geodesic transversal problem finds a minimum vertex set that hits all possible maximal geodesics, providing a more focused framework for applications such as load balancing, sensor placement, and propagation control in large scale network, thus making it a more powerful model for addressing challenges in modern network analysis. We begin by quoting a few preliminary results in section 2, then deduce a lower bound for the geodesic transversal number of an arbitrary graph $G$  in section 3 and compute the gt-number of Butterfly and Hypercube networks in sections 4 and 5 respectively. 

 \section{Preliminaries}
 The propositions stated below were established in \cite{geo.cover.butterfly} while investigating the geodesic cover problem of $BF(r)$. We apply these propositions in section $4$, to investigate the geodesic transversal number of $BF(r)$.

\begin{prop}
    \label{samelevel}\normalfont{\cite{geo.cover.butterfly}}
    \em{A geodesic of $BF(r)$ contains at most two vertices of level $0$ and at most two  vertices of level $r$. Moreover, if a geodesic contains two vertices of level 0, then they are the ends of the geodesic (and similarly for level r)}. 
\end{prop}
    
\begin{prop}
    \label{sameleve.max}\normalfont{\cite{geo.cover.butterfly}}
    \em{If both end vertices of a geodesic $P$ of $BF(r)$ are either at level $0$ or at level $r$, then $P$ is maximal}.
\end{prop}

\begin{prop}
    \normalfont{\cite{geo.cover.butterfly}}
    \em{If a geodesic $P$ of $BF(r)$ covers three vertices of degree $2$, then the end vertices of $P$ are of degree $2$, and $P$ is maximal}.
\end{prop}
\section{Lower bound}
 In this section, we study the behaviour of the maximal geodesics of induced subgraphs in its parent graph. We first consider the case, where a maximal geodesic of $G$ contains two or more maximal geodesics of an induced subgraph or distinct induced subgraphs refer Figures \ref{gen.example.1} and \ref{gen.example.2}. In this case, we prove that the gt-number of the parent graph $G$ is bounded above by the sum of the gt-number of its induced subgraphs (Lemma \ref{upper}). Then we prove the most significant result of this paper, called the Geodesic Transversal Lemma in which the gt-number of $G$ is bounded below by the sum of the gt-numbers of its induced subgraphs when no maximal geodesic of $G$ contains two or more maximal geodesics of same induced subgraph or distinct induced subgraphs of $G$. 

 \begin{lem}
     \label{upper}
    Let $G(V,E)$ be a connected graph and $G[V_1], G[V_2],\ldots,G[V_k]$ be induced subgraphs of $G$ such that $V_i\cap V_j= \emptyset$, $1\leq i, j \leq k, i\neq j$. Let $S_i$ be the set of all possible maximal geodesics of $G[V_i]$ respectively and $S=S_1\cup S_2\cup \ldots \cup S_k$. If there exist no maximal geodesics of $G$ that lies entirely in $G[V\backslash(V_1\cup V_2\cup\ldots\cup V_k)]$ then, $gt(G)\leq gt(G[V_1])+gt(G[V_2])+\ldots +gt(G[V_k])$.
 \end{lem}
 \begin{proof}
     Let $T_1, T_2,\ldots ,T_k$ be the gt-sets of $G[V_1], G[V_2],\ldots,G[V_k]$ respectively, $T=T_1\cup T_2\cup \ldots \cup T_k$ and let $S'$ be the set of all maximal geodesics of $G$ that lie entirely in $G[V\backslash (V_1\cup \ldots\cup V_n)]$. By our hypothesis, $S'=\emptyset$. Then the maximal geodesics of $G$ are either maximal geodesics of $G[V_i],1\leq i\leq k$ or includes two or more maximal geodesics of $S$. Let $P$ be a maximal geodesic of $G$ that contains maximal geodesics say $P_1, P_2,\ldots ,P_l\subseteq S$, then the set of vertices say $\{v_1, v_2,\ldots ,v_l\}$ in $T$ that hits $P_1, P_2,\ldots ,P_l$ respectively also hits the maximal geodesic $P$. Let $V'=\{v_1, v_2,\ldots,v_l\}$. Without loss of generality, assume that the remaining maximal geodesics in $S \backslash \{P_1, P_2,\ldots,P_l\}$ will either be maximal geodesic or extended maximal geodesic of $G$. Depending on the number of maximal geodesics hit by each vertex $v_i$ of $V'$ , we consider the following cases: \\
    \textit{Case i}: If each $v_i\in V'$ hits exactly one geodesic, then $T-\{v_2,v_3,\ldots,v_l\}$ is a gt-set of $G$ and hence $gt(G)< gt(G[V_1])+gt(G[V_2])+\ldots +gt(G[V_k])$.\\
    \textit{Case ii}: For $q<l$, if every vertex of $\{v_{i_1},v_{i_2},\ldots,v_{i_q}\}\subset V' $ hits exactly one maximal geodesic of $S$, then $T-\{v_{i_1},v_{i_2},\ldots,v_{i_q}\}$ is a gt-set of $G$ and hence $gt(G)< gt(G[V_1])+gt(G[V_2])+\ldots+gt(G[V_k])$.\\
    \textit{Case iii}: Suppose each vertex in $V'$ hits at least two maximal geodesics in $S$. Without loss of generality, consider $v_i$, $1\leq i\leq l$ that hits $P_i \in \{P_1,P_2,...,P_l\}$ and $P'\in S \backslash \{P_1, P_2,\ldots,P_l\}$. Then $T-v_i$ is not a geodesic transversal of $G$ since $P'$ does not have any hit vertex. Thus, $T$ is minimal and we have $gt(G)= gt(G[V_1])+gt(G[V_2])+\ldots+gt(G[V_k])$. Hence the proof. 
 \end{proof}
 
\begin{figure}[htbp]
    \centering
\begin{tikzpicture}[scale=0.85,
    x=0.01cm,
    y=-0.01cm,
    vertex/.append style={minimum size=2mm},
    redvertex/.append style={minimum size=2mm},
    bluevertex/.append style={minimum size=2mm}
]

\begin{scope}[shift={(0,0)}]

\node[vertex]    (a1) at (116,21) {};
\node[redvertex] (a2) at (209,113) {};
\node[vertex]    (a3) at (250,141) {};
\node[vertex]    (a4) at (293,170) {};
\node[vertex]    (a5) at (339,199) {};
\node[vertex]    (a6) at (339,114) {};
\node[vertex]    (a7) at (433,163) {};
\node[vertex]    (a8) at (381,227) {};
\node[vertex]    (a9) at (423,255) {};
\node[vertex]    (a10) at (462,281) {};
\node[vertex]    (a11) at (499,304) {};
\node[redvertex]    (a12) at (534,328) {};
\node[vertex]    (a13) at (638,328) {};
\node[vertex]    (a14) at (534,425) {};
\node[vertex]    (a15) at (413,488) {};
\node[vertex]    (a16) at (203,453) {};
\node[vertex]    (a17) at (204,622) {};
\node[vertex]    (a18) at (294,660) {};
\node[vertex]    (a19) at (317,576) {};
\node[vertex]    (a20) at (352,753) {};
\node[vertex]    (a21) at (409,702) {};
\node[vertex]    (a22) at (412,811) {};
\node[vertex]    (a23) at (79,170) {};
\node[vertex]    (a24) at (59,308) {};

\draw[reddash] (a1)--(a2);
\draw[reddash] (a2)--(a23);
\draw[reddash] (a23)--(a24);
\draw[reddash] (a24)--(a16);
\draw[reddash] (a16)--(a15);
\draw[reddash] (a15)--(a14);
\draw[reddash] (a14)--(a12);
\draw[edge] (a16)--(a17)--(a18);
\draw[edge] (a18)--(a19);
\draw[edge] (a18)--(a20);
\draw[edge] (a20)--(a21);
\draw[edge] (a20)--(a22);
\draw[edge] (a2)--(a3)--(a4)--(a5);
\draw[edge] (a5)--(a6);
\draw[edge] (a5)--(a7);
\draw[edge] (a5)--(a8)--(a9)--(a10)--(a11)--(a12);
\draw[reddash] (a12)--(a13);

\node at (119,96) {$P_1$};
\node at (590,381) {$P_2$};
\node at (275,888) {$G$};

\end{scope}

\begin{scope}[shift={(100,0)}]

\node[vertex]    (b1) at (838,183) {};
\node[redvertex] (b2) at (893,293) {};
\node[vertex]    (b3) at (761,304) {};
\node[vertex]    (b4) at (922,330) {};
\node[vertex]    (b5) at (951,368) {};
\node[vertex]    (b6) at (985,407) {};
\node[vertex]    (b7) at (1009,330) {};
\node[vertex]    (b8) at (1079,403) {};
\node[vertex]    (b9) at (1014,445) {};
\node[vertex]    (b10) at (1043,482) {};
\node[vertex]    (b11) at (1070,516) {};
\node[vertex]    (b12) at (1096,547) {};
\node[redvertex] (b13) at (1121,579) {};
\node[vertex]    (b14) at (1217,611) {};
\node[vertex]    (b15) at (1093,668) {};

\draw[rededge] (b1)--(b2);
\draw[rededge] (b3)--(b2);
\draw[edge] (b2)--(b4)--(b5)--(b6);
\draw[edge] (b6)--(b7);
\draw[edge] (b6)--(b8);
\draw[edge] (b6)--(b9)--(b10)--(b11)--(b12)--(b13);
\draw[rededge] (b13)--(b14);
\draw[rededge] (b13)--(b15);

\node at (814,254) {$P_1$};
\node at (1165,650) {$P_2$};
\node at (949,888) {$G\,[V_1]$};

\end{scope}

\begin{scope}[shift={(200,0)}]

\node[vertex] (c1) at (1312,315) {};
\node[vertex] (c2) at (1401,353) {};
\node[vertex] (c3) at (1424,268) {};
\node[vertex] (c4) at (1459,446) {};
\node[vertex] (c5) at (1517,395) {};
\node[vertex] (c6) at (1518,503) {};

\draw[edge] (c1)--(c2);
\draw[edge] (c2)--(c3);
\draw[edge] (c2)--(c4);
\draw[edge] (c4)--(c5);
\draw[edge] (c4)--(c6);

\node at (1447,888) {$G\,[V_2]$};

\end{scope}

\end{tikzpicture}
\caption{Maximal geodesic $P$ (dotted lines) in $G$ that includes two maximal geodesics $P_1$ and $P_2$ of $G[V_1]$.}
\label{gen.example.1}
\end{figure}
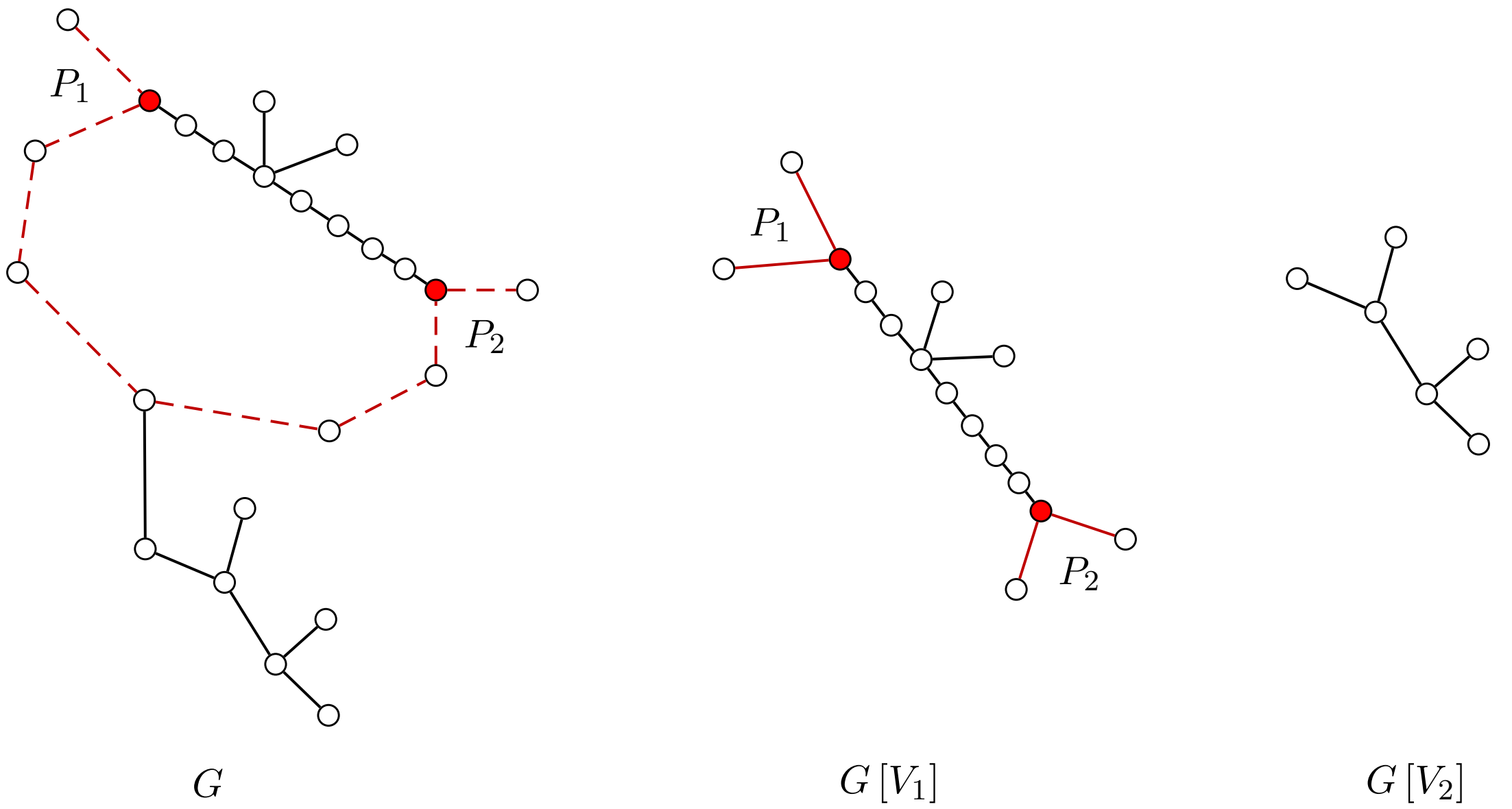

\begin{figure}[htbp]
    \centering
\begin{tikzpicture}[scale=0.85,
    x=0.01cm,
    y=-0.01cm,
    vertex/.append style={minimum size=2mm},
    redvertex/.append style={minimum size=2mm},
    bluevertex/.append style={minimum size=2mm}
]

\begin{scope}[shift={(0,0)}]

\node[vertex]    (d1)  at (184,74) {};
\node[redvertex] (d2)  at (268,158) {};
\node[vertex]    (d3)  at (306,183) {};
\node[vertex]    (d4)  at (343,208) {};
\node[vertex]    (d5)  at (386,235) {};
\node[vertex]    (d6)  at (386,158) {};
\node[vertex]    (d7)  at (470,203) {};
\node[vertex]    (d8)  at (423,260) {};
\node[vertex]    (d9)  at (461,286) {};
\node[vertex]    (d10) at (497,308) {};
\node[vertex]    (d11) at (530,329) {};
\node[vertex]    (d12) at (561,351) {};
\node[vertex]    (d13) at (656,351) {};
\node[vertex]    (d14) at (561,439) {};
\node[vertex]    (d15) at (452,495) {};
\node[vertex]    (d16) at (263,464) {};
\node[vertex]    (d17) at (263,617) {};
\node[bluevertex] (d18) at (345,650) {};
\node[vertex]    (d19) at (366,575) {};
\node[vertex]    (d20) at (397,734) {};
\node[vertex]    (d21) at (449,689) {};
\node[vertex]    (d22) at (452,787) {};
\node[vertex]    (d23) at (150,208) {};
\node[vertex]    (d24) at (133,333) {};

\draw[reddash] (d1)--(d2);
\draw[reddash] (d2)--(d23);

\draw[blackdash] (d23)--(d24);
\draw[blackdash] (d24)--(d16);
\draw[blackdash] (d16)--(d17);

\draw[bluedash] (d17)--(d18);
\draw[bluedash] (d18)--(d19);

\draw[edge] (d2)--(d3)--(d4)--(d5);
\draw[edge] (d5)--(d6);
\draw[edge] (d5)--(d7);
\draw[edge] (d5)--(d8)--(d9)--(d10)--(d11)--(d12);
\draw[edge] (d12)--(d13);
\draw[edge] (d12)--(d14)--(d15);
\draw[edge] (d16)--(d15);
\draw[edge] (d18)--(d20);
\draw[edge] (d20)--(d21);
\draw[edge] (d20)--(d22);

\node at (195,140) {$P_1$};
\node at (320,602) {$P_2$};
\node at (298,866) {$G$};

\end{scope}

\begin{scope}[shift={(100,0)}]

\node[vertex]    (e1)  at (835,220) {};
\node[redvertex] (e2)  at (885,319) {};
\node[vertex]    (e3)  at (766,329) {};
\node[vertex]    (e4)  at (912,353) {};
\node[vertex]    (e5)  at (938,387) {};
\node[vertex]    (e6)  at (968,423) {};
\node[vertex]    (e7)  at (990,354) {};
\node[vertex]    (e8)  at (1053,419) {};
\node[vertex]    (e9)  at (994,456) {};
\node[vertex]    (e10) at (1021,490) {};
\node[vertex]    (e11) at (1045,521) {};
\node[vertex]    (e12) at (1070,550) {};
\node[vertex]    (e13) at (1090,577) {};
\node[vertex]    (e14) at (1178,607) {};
\node[vertex]    (e15) at (1065,657) {};

\draw[rededge] (e1)--(e2);
\draw[rededge] (e3)--(e2);
\draw[edge] (e2)--(e4)--(e5)--(e6);
\draw[edge] (e6)--(e7);
\draw[edge] (e6)--(e8);
\draw[edge] (e6)--(e9)--(e10)--(e11)--(e12)--(e13);
\draw[edge] (e13)--(e14);
\draw[edge] (e13)--(e15);

\node at (814,282) {$P_1$};
\node at (967,866) {$G\,[V_1]$};

\end{scope}

\begin{scope}[shift={(200,0)}]

\node[vertex]     (f1) at (1263,339) {};
\node[bluevertex] (f2) at (1344,373) {};
\node[vertex]     (f3) at (1364,297) {};
\node[vertex]     (f4) at (1397,457) {};
\node[vertex]     (f5) at (1447,411) {};
\node[vertex]     (f6) at (1450,509) {};

\draw[blueedge] (f1)--(f2);
\draw[blueedge] (f2)--(f3);
\draw[edge] (f2)--(f4);
\draw[edge] (f4)--(f5);
\draw[edge] (f4)--(f6);

\node at (1311,319) {$P_2$};
\node at (1411,866) {$G\,[V_2]$};

\end{scope}

\end{tikzpicture}

\caption{Maximal geodesic $P$ (dotted lines) in $G$ containing two maximal geodesics $P_1$ of $G[V_1]$ and $P_2$ of $G[V_2]$.}
\label{gen.example.2}
\end{figure}
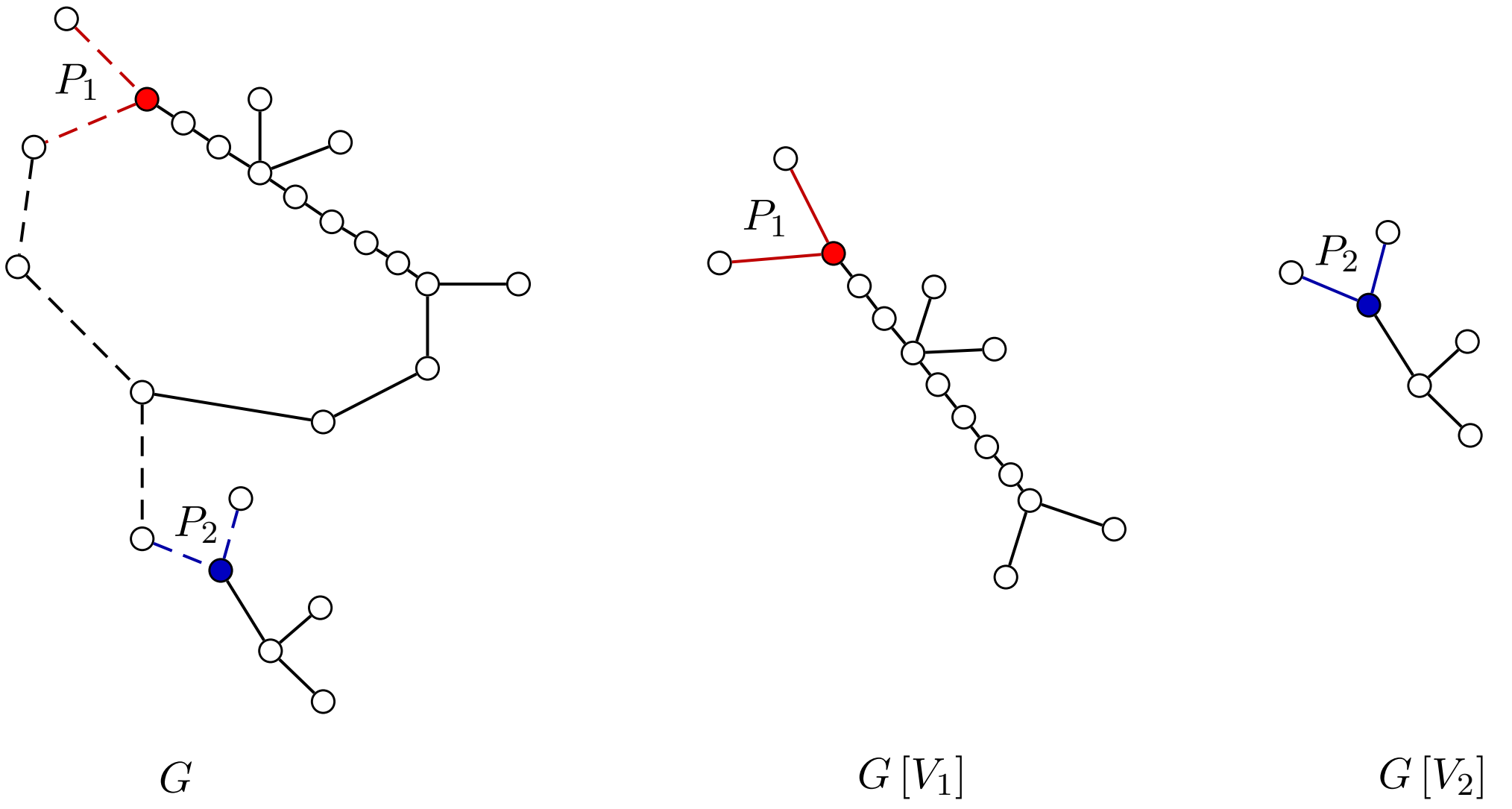

\begin{lem}
    \label{GTL}\normalfont{(Geodesic Transversal lemma)}.
    \em{Let $G(V,E)$ be a connected graph and its induced subgraphs be denoted as $G[V_1], G[V_2],\ldots,G[V_k]$ such that $V_i\cap V_j= \emptyset$, $1\leq i, j \leq k, i\neq j$. Let $S_i$ be the set of all possible maximal geodesics of $G[V_i], 1\leq i\leq k$, respectively, and let $S=S_1\cup S_2\cup \ldots\cup S_k$. If no maximal geodesic $P$ in $G$ includes two maximal geodesics of $S$ as its sub-path, then $gt(G)\geq gt(G[V_1])+gt(G[V_2])+\ldots+gt(G[V_k])$}.
\end{lem}
\begin{proof}
    Let $T_i$ be the gt-set of $G[V_i], 1\leq i\leq k$ and $T=T_1\cup T_2\cup \ldots \cup T_k$. Consider the induced subgraph $G[V\backslash (V_1\cup V_2\cup\ldots\cup V_k)]$ and let $T'$ be its gt-set. Denote by $S'$ the set of all maximal geodesics of $G$ that lie entirely in $G[V\backslash (V_1\cup V_2\cup\ldots\cup V_k)]$. By our hypothesis, since no maximal geodesic in $G$ includes two maximal geodesics of $S$ as its sub-path, every maximal geodesic in $S$ will either be a maximal geodesic or sub-path of a maximal geodesic in  $G$. The set $S'$ may be empty or non empty, depending on which we have the following cases: \\
\textit{Case i}: Every maximal geodesic $P_i \in S$ is a maximal geodesic in $G$. In this case, 
\begin{enumerate}
    \item [a.] If $S'=\emptyset$, then $T$ will be a gt-set of $G$ and $gt(G)=gt(G[V_1])+gt(G[V_2])+\ldots + gt(G[V_k])$.
    \item [b.] If $S'\neq \emptyset$, then $T\cup T'$ is a gt-set of $G$ and $gt(G)> gt(G[V_1])+gt(G[V_2])+\ldots+ gt(G[V_k])$.
\end{enumerate}
\textit{Case ii}: There exists at least one maximal geodesic, say $P_i\in S$ which is not maximal in $G$. Clearly, the maximal geodesic $P_i$ will be a sub-path of a maximal geodesic, say $P$ in $G$. Then the vertex, $v_i\in T$ that hits the maximal geodesic $P_i$ will also hit the maximal geodesic $P$. Further, 
\begin{enumerate}
    \item [a.] If $S'=\emptyset,$ then $T$ will be a gt-set of $G$ and $gt(G)=gt(G[V_1])+gt(G[V_2])+\ldots +gt(G[V_k])$.
    \item [b.] If $S'\neq \emptyset$,then $T\cup T'$ is a gt-set of $G$ and $gt(G)>gt(G[V_1])+gt(G[V_2])+\ldots+gt(G[V_k])$.
\end{enumerate}
In cases i(a) and ii(a), since there exists no maximal geodesic of $G$ that lies entirely in $G[V\backslash (V_1\cup V_2\cup\ldots\cup V_k)]$, the maximal geodesics of the induced subgraphs are the maximal geodesics or extended maximal geodesics of $G$ and hence $gt(G)=gt(G[V_1])+gt(G[V_2])+\ldots +gt(G[V_k])$, while in cases i(b) and ii(b) there exists some maximal geodesics that lie entirely in $G[V\backslash (V_1\cup V_2\cup\ldots\cup V_k)]$, and also the maximal geodesics of the induced subgraphs remain as maximal geodesics or extended maximal geodesics of $G$, which implies $gt(G)>gt(G[V_1])+gt(G[V_2])+\ldots+gt(G[V_k])$. Combining these results, we have $gt(G)\geq gt(G[V_1])+gt(G[V_2])+\ldots+gt(G[V_k])$.
\end{proof}

\begin{cor}
    \label{gtlc}
	Let $G=(V,E)$ be a connected graph, $G[V_1]$ and $G[V_2]$ be induced subgraphs of $G$ such that $V_1\cap V_2= \emptyset$ and $S_1$, $S_2$ be the set of all possible maximal geodesics of $G[V_1]$ and $G[V_2]$ respectively, with $S=S_1\cup S_2$. If no maximal geodesic $P$ in $G$ includes two maximal geodesics of $S$ as its sub-path, then $gt(G)\geq gt(G[V_1])+gt(G[V_2])$.
\end{cor}
\begin{proof}
    The proof follows directly from the Geodesic Transversal lemma for $k=2$. 
\end{proof}
\section{Butterfly networks}\label{BF(r)}

In this section, we determine the geodesic transversal number of Butterfly network, a class of recursively defined graphs with wide applications in parallel computing and communication networks using the Geodesic Transversal lemma.

The butterfly network of $r-$dimension is denoted as $BF(r)$. It is a $(r+1)-$partite graph. The vertices of $BF(r)$ correspond to a set of pairs $[s,l]$, where $l$ denotes the level of the vertex, $0 \leq l \leq r$ and $s$ is a $r-$bit binary number. Two vertices $[s,l]$ and $[s', l']$ are adjacent in $BF(r)$ if and only if $l'=l+1$ and either $s$ and $s'$ are identical or $s$ and $s'$ differ in the $l$-th bit~\cite{iso.path.part.butterfly}. $BF(r)$ contains two copies of $BF(r-1)$. Throughout this section, we denote the two copies of $BF(r-1)$ as $BF'(r-1)$ and $BF''(r-1)$. The following Lemmas~\ref{case2},~\ref{case3} and~\ref{case1}, describe the maximal geodesics of $BF(r-1)$ in $BF(r)$. 

\begin{figure}[htbp]
\centering
\begin{tikzpicture}[scale=0.85,
    x=0.01cm,
    y=-0.01cm,
    vertex/.append style={minimum size=2mm},
    redvertex/.append style={minimum size=2mm},
    bluevertex/.append style={minimum size=2mm}
]

\draw[black, dashed, line width=0.8pt]
    (214,139) rectangle (734,686);

\draw[black, dashed, line width=0.8pt]
    (760,139) rectangle (1281,686);

\node[vertex] (a00) at (269,194) {};
\node[vertex] (a01) at (406,194) {};
\node[vertex] (a02) at (542,194) {};
\node[vertex] (a03) at (679,194) {};
\node[vertex] (a04) at (815,194) {};
\node[vertex] (a05) at (952,194) {};
\node[vertex] (a06) at (1089,194) {};
\node[vertex] (a07) at (1225,194) {};

\node[vertex] (b00) at (269,413) {};
\node[vertex] (b01) at (406,413) {};
\node[vertex] (b02) at (542,413) {};
\node[vertex] (b03) at (679,413) {};
\node[vertex] (b04) at (815,413) {};
\node[vertex] (b05) at (952,413) {};
\node[vertex] (b06) at (1089,413) {};
\node[vertex] (b07) at (1225,413) {};

\node[vertex] (c00) at (269,631) {};
\node[vertex] (c01) at (406,631) {};
\node[vertex] (c02) at (542,631) {};
\node[vertex] (c03) at (679,631) {};
\node[vertex] (c04) at (815,631) {};
\node[vertex] (c05) at (952,631) {};
\node[vertex] (c06) at (1089,631) {};
\node[vertex] (c07) at (1225,631) {};

\node[vertex] (d00) at (269,904) {};
\node[vertex] (d01) at (406,904) {};
\node[vertex] (d02) at (542,904) {};
\node[vertex] (d03) at (679,904) {};
\node[vertex] (d04) at (815,904) {};
\node[vertex] (d05) at (952,904) {};
\node[vertex] (d06) at (1089,904) {};
\node[vertex] (d07) at (1225,904) {};

\draw[edge] (a00)--(b00);
\draw[edge] (a01)--(b01);
\draw[edge] (a02)--(b02);
\draw[edge] (a03)--(b03);
\draw[edge] (a04)--(b04);
\draw[edge] (a05)--(b05);
\draw[edge] (a06)--(b06);
\draw[edge] (a07)--(b07);

\draw[edge] (a00)--(b01);
\draw[edge] (a01)--(b00);
\draw[edge] (a02)--(b03);
\draw[edge] (a03)--(b02);
\draw[edge] (a04)--(b05);
\draw[edge] (a05)--(b04);
\draw[edge] (a06)--(b07);
\draw[edge] (a07)--(b06);

\draw[edge] (b00)--(c00);
\draw[edge] (b01)--(c01);
\draw[edge] (b02)--(c02);

\draw[edge] (b04)--(c04);
\draw[edge] (b05)--(c05);
\draw[edge] (b06)--(c06);
\draw[edge] (b07)--(c07);

\draw[edge] (b00)--(c02);
\draw[edge] (b01)--(c03);
\draw[edge] (b02)--(c00);

\draw[edge] (b04)--(c06);
\draw[edge] (b05)--(c07);
\draw[edge] (b06)--(c04);
\draw[edge] (b07)--(c05);

\draw[edge] (c00)--(d00);

\draw[edge] (c02)--(d02);

\draw[edge] (c04)--(d04);
\draw[edge] (c05)--(d05);
\draw[edge] (c06)--(d06);
\draw[edge] (c07)--(d07);

\draw[edge] (c00)--(d04);
\draw[edge] (c01)--(d05);
\draw[edge] (c02)--(d06);
\draw[edge] (c03)--(d07);
\draw[edge] (c04)--(d00);
\draw[edge] (c05)--(d01);
\draw[edge] (c06)--(d02);
\draw[edge] (c07)--(d03);


\draw[rededge] (c01)--(b03);
\draw[rededge] (b03)--(c03);

\draw[reddash] (c01)--(d01);
\draw[reddash] (c03)--(d03);

\node[anchor=east] at (398,640) {$u$};
\node[anchor=east] at (669,640) {$v$};

\node[anchor=east] at (398,915) {$u'$};
\node[anchor=east] at (669,915) {$v'$};

\node[anchor=east] at (178,194) {$\ell=0$};
\node[anchor=east] at (178,413) {$\ell=1$};
\node[anchor=east] at (178,631) {$\ell=2$};
\node[anchor=east] at (178,904) {$\ell=3$};

\node at (269,968) {$000$};
\node at (406,968) {$001$};
\node at (542,968) {$010$};
\node at (679,968) {$011$};
\node at (815,968) {$100$};
\node at (952,968) {$101$};
\node at (1089,968) {$110$};
\node at (1225,968) {$111$};

\end{tikzpicture}
    \caption{Maximal geodesic $P(u,v)$ of $BF(r-1)$ extended to $u'$ and $v'$ respectively} \label{fig3:4.1case1&2}
\end{figure}
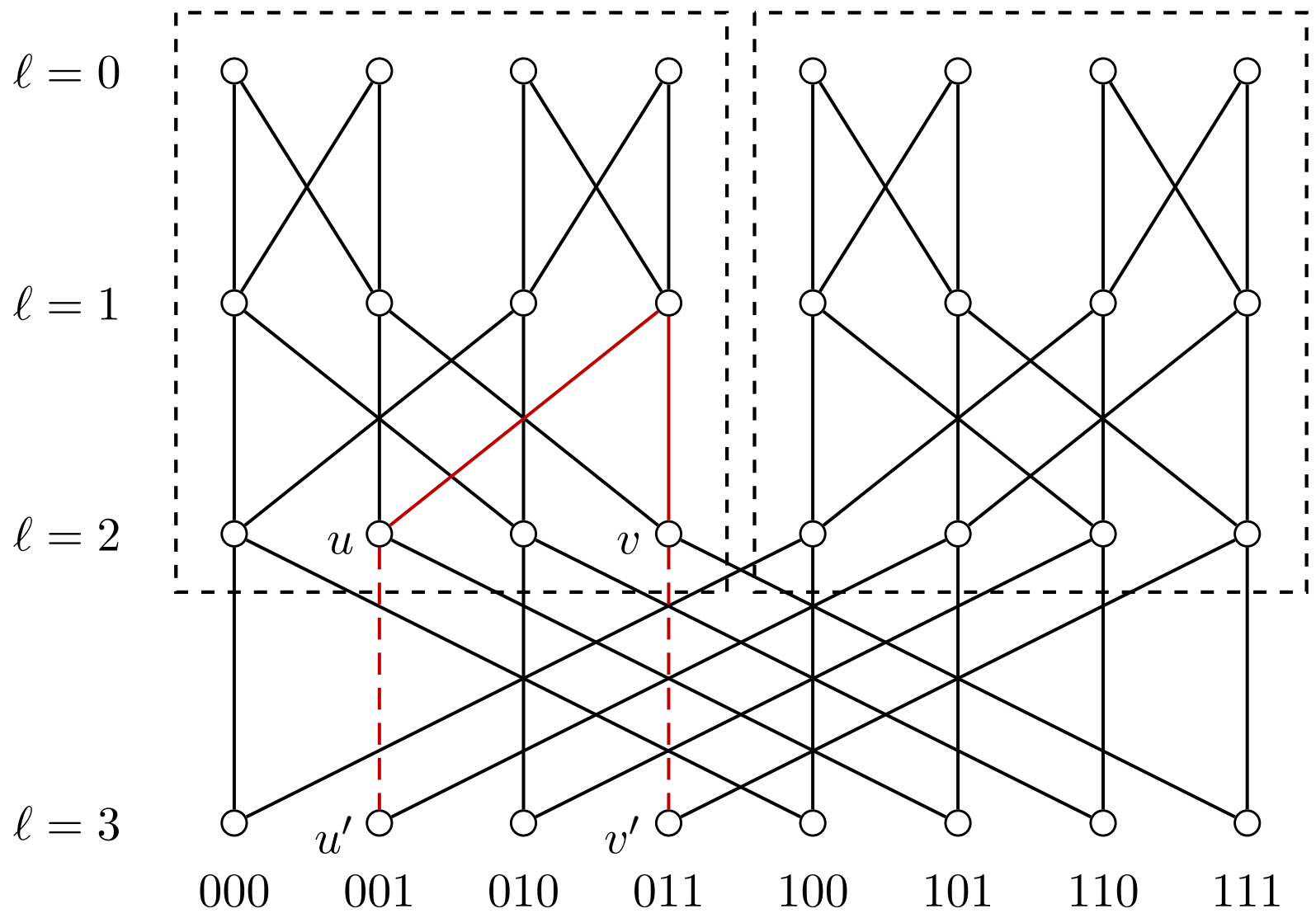  

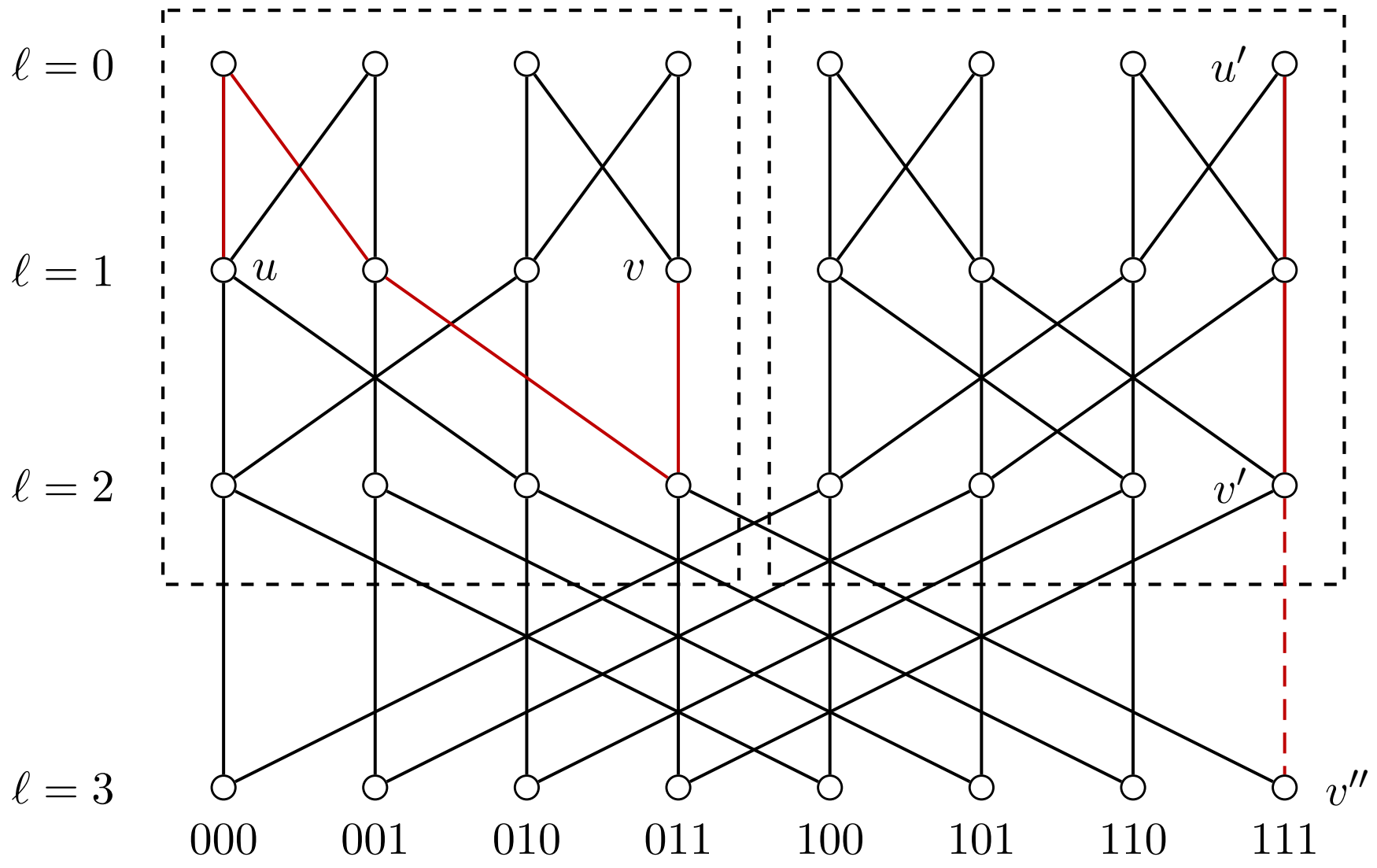
\begin{figure}[htbp]
\centering
\begin{tikzpicture}[scale=0.85,
    x=0.01cm,
    y=-0.01cm,
    vertex/.append style={minimum size=2mm},
    redvertex/.append style={minimum size=2mm},
    bluevertex/.append style={minimum size=2mm}
]

\draw[black, dashed, line width=0.8pt]
    (198,116) rectangle (768,684);

\draw[black, dashed, line width=0.8pt]
    (798,116) rectangle (1367,684);

\node[vertex] (a0) at (258,169) {};
\node[vertex] (a1) at (408,169) {};
\node[vertex] (a2) at (558,169) {};
\node[vertex] (a3) at (708,169) {};
\node[vertex] (a4) at (858,169) {};
\node[vertex] (a5) at (1008,169) {};
\node[vertex] (a6) at (1158,169) {};
\node[vertex] (a7) at (1308,169) {};

\node[vertex] (b0) at (258,373) {};
\node[vertex] (b1) at (408,373) {};
\node[vertex] (b2) at (558,373) {};
\node[vertex] (b3) at (708,373) {};
\node[vertex] (b4) at (858,373) {};
\node[vertex] (b5) at (1008,373) {};
\node[vertex] (b6) at (1158,373) {};
\node[vertex] (b7) at (1308,373) {};

\node[vertex] (c0) at (258,586) {};
\node[vertex] (c1) at (408,586) {};
\node[vertex] (c2) at (558,586) {};
\node[vertex] (c3) at (708,586) {};
\node[vertex] (c4) at (858,586) {};
\node[vertex] (c5) at (1008,586) {};
\node[vertex] (c6) at (1158,586) {};
\node[vertex] (c7) at (1308,586) {};

\node[vertex] (d0) at (258,885) {};
\node[vertex] (d1) at (408,885) {};
\node[vertex] (d2) at (558,885) {};
\node[vertex] (d3) at (708,885) {};
\node[vertex] (d4) at (858,885) {};
\node[vertex] (d5) at (1008,885) {};
\node[vertex] (d6) at (1158,885) {};
\node[vertex] (d7) at (1308,885) {};


\draw[edge] (a0)--(b0);
\draw[edge] (a1)--(b1);
\draw[edge] (a2)--(b2);
\draw[edge] (a3)--(b3);
\draw[edge] (a4)--(b4);
\draw[edge] (a5)--(b5);
\draw[edge] (a6)--(b6);
\draw[edge] (a7)--(b7);

\draw[rededge] (a0)--(b1);
\draw[edge]    (a1)--(b0);

\draw[edge] (a2)--(b3);
\draw[edge] (a3)--(b2);

\draw[edge] (a4)--(b5);
\draw[edge] (a5)--(b4);

\draw[edge] (a6)--(b7);
\draw[edge] (a7)--(b6);


\draw[edge] (b0)--(c0);
\draw[edge] (b1)--(c1);
\draw[edge] (b2)--(c2);

\draw[edge] (b4)--(c4);
\draw[edge] (b5)--(c5);
\draw[edge] (b6)--(c6);
\draw[rededge] (b7)--(c7);

\draw[edge] (b0)--(c2);
\draw[edge] (b2)--(c0);

\draw[edge] (b4)--(c6);
\draw[edge] (b5)--(c7);
\draw[edge] (b6)--(c4);
\draw[edge] (b7)--(c5);


\draw[edge] (c0)--(d0);
\draw[edge] (c1)--(d1);
\draw[edge] (c2)--(d2);
\draw[edge] (c3)--(d3);
\draw[edge] (c4)--(d4);
\draw[edge] (c5)--(d5);
\draw[edge] (c6)--(d6);
\draw[reddash] (c7)--(d7);

\draw[edge] (c0)--(d4);
\draw[edge] (c1)--(d5);
\draw[edge] (c2)--(d6);
\draw[edge] (c3)--(d7);

\draw[edge] (c4)--(d0);
\draw[edge] (c5)--(d1);
\draw[edge] (c6)--(d2);
\draw[edge] (c7)--(d3);


\draw[rededge] (a0)--(b0);

\draw[rededge] (b1)--(c3);
\draw[rededge] (b3)--(c3);

\draw[rededge] (a7)--(b7);
\draw[rededge] (b7)--(c7);

\node[anchor=east] at (166,169) {$\ell=0$};
\node[anchor=east] at (166,373) {$\ell=1$};
\node[anchor=east] at (166,586) {$\ell=2$};
\node[anchor=east] at (166,885) {$\ell=3$};

\node[anchor=west] at (270,373) {$u$};
\node[anchor=east] at (692,373) {$v$};

\node[anchor=east] at (1288,169) {$u'$};
\node[anchor=east] at (1288,586) {$v'$};
\node[anchor=west] at (1332,885) {$v''$};

\node at (258,936) {$000$};
\node at (408,936) {$001$};
\node at (558,936) {$010$};
\node at (708,936) {$011$};
\node at (858,936) {$100$};
\node at (1008,936) {$101$};
\node at (1158,936) {$110$};
\node at (1308,936) {$111$};

\end{tikzpicture}
    \caption{A maximal geodesic $P=(u,v)$ where $d(u)=d(v)=4$ and a geodesic $Q=(u',v')$ where $d(u')=2$ and $d(v')=4$}
    \label{fig4:4.2&4.3}
    \end{figure}

\begin{lem}
    \label{case2}
   A maximal geodesic $P$ in $BF(r-1)$ whose both end vertices are of degree two is either a maximal geodesic or an extended maximal geodesic of $BF(r)$.
\end{lem}
\begin{proof}
    Let $u$ and $v$ be the end vertices of $P$, a maximal geodesic in $BF(r-1)$ where $d(u)=d(v)=2$. Clearly, both $u$ and $v$ either belong to the level $0$ or the level $r-1$ in $BF(r-1)$. For otherwise, if $P(u,v)$ is a maximal geodesic with $u=(s,0)$ and $v=(s',r-1)/(s,r-1)$ then there exists a path $P'=P\cup ((s',r-1),(s',r-2))/P\cup ((s,r-1),(s',r-2))$ respectively which is a geodesic of $BF(r-1)$ containing $P$, contradicting the fact $P$ is maximal.  Hence we proceed by considering the following two cases:\\
\noindent\textit{Case i} : $u=(s,0)$ and $v=(s',0)$\\
By Proposition~\ref{sameleve.max}, $P$ is a maximal geodesic of $BF(r)$. \\
\noindent\textit{Case ii} : $u=(s,r-1)$ and $v=(s',r-1)$\\
The geodesic $P$ can be extended to vertices of degree two of $BF(r)$, namely, $u'=(s,r)$ and $v'=(s',r)$, which is a maximal geodesic in $BF(r)$, refer Figure~\ref{fig3:4.1case1&2}.
\end{proof}

\begin{lem}
    \label{case3}
    A maximal geodesic $P$ of $BF(r-1)$ whose both end vertices are of degree four is a maximal geodesic in $BF(r)$.
\end{lem}
\begin{proof}
    Assume the contrary. If $P$ is not maximal in $BF(r)$, then there exists a maximal geodesic $P'$ in $BF(r)$ that includes $P$. Now, $P'\cap BF(r-1)$ is a geodesic in $BF(r-1)$ containing $P$ which contradicts the fact that $P$ is a maximal in $BF(r-1)$, refer Figure~\ref{fig4:4.2&4.3}. 
\end{proof}

\begin{lem}
    \label{case1}
 If $P(u',v')$ is a geodesic in $BF(r)$ such that $d(u')=2$ and $d(v')=4$, then $P$ is not maximal in $BF(r)$.
\end{lem}

    \begin{proof}
    We may assume without loss of generality that $u'=(s,0)$ and $v'=(s,j)$, $1\leq j \leq r-1$. Then $P(u',v')$ can be extended to $P(u',v'')$, $v''=(s,j+1)$, which proves that $P$ is not maximal, refer Figure~\ref{fig4:4.2&4.3}.
   
\end{proof}

\begin{thm}
    \label{BF-gt}
    If $BF(r)$ is a butterfly network of dimension $r$, then $gt(BF(r))=2^{r+1},r\geq 3$.
\end{thm}
\begin{proof}
    Since every maximal geodesic of $BF(r)$ includes at least one vertex of degree two, the set of vertices at level $0$ and level $r$ is a $gt-set$ of $BF(r)$. Hence, $gt(BF(r))\leq 2^{r+1}$. \\
    The proof of the other inequality is by induction on $r$. Let $r=3$ and $S$ be a geodesic transversal of $BF(3)$. Depending on the vertices of $S$ we have the following cases. \\
\textit{Case i} : $S$ contains only vertices of degree two. Suppose that $(000,0)\notin S$, then the maximal geodesic $P:(001,1)\rightarrow (000,0)\rightarrow (000,1)\rightarrow (000,2)\rightarrow (010,1)$ does not have any hit vertex in $S$ and hence $gt(BF(3))\geq 2^3+2^3=16$. \\
\textit{Case ii} : $S$ contains only vertices of degree four. Suppose that $(000,1)\notin S$, then the maximal geodesic $P: (000,0)\rightarrow(000,1)\rightarrow(001,0)$ does not have any hit vertex in $S$ and hence $gt(BF(3))\geq 2^3+2^3=16$.\\
\textit{Case iii} : $S$ contains vertices of degree two and four. Let $S$ be the set of all vertices of level $0$ and level $2$ such that all the maximal geodesics of length two are hit. Suppose $(111,2)\notin S$ then the maximal geodesic $P:(111,3)\rightarrow (111,2)\rightarrow (011,3)$ does not have any hit vertex in $S$ and hence $gt(BF(3))\geq 2^3+2^3=16$. Thus, the base case for $r=3$ is verified.

Assume the result is true for $r=k-1$, that is $gt(BF(k-1))\geq 2^{(k-1)+1}$. We claim that $gt(BF(k))\geq 2^{k+1}$. Let $G[V_1]$ and $G[V_2]$ be the induced subgraphs $BF(k-1)'$ and $BF(k-1)''$ and let $S_1$ and $S_2$ denote the maximal geodesics in $G[V_1]$ and $G[V_2]$ respectively. By Lemmas \ref{case2}, \ref{case3} and \ref{case1} we observe that no maximal geodesic of $BF(k)$ contains two maximal geodesics either from $S_1$ (respectively $S_2$) or one from each sets $S_1$ and $S_2$. Hence, by applying Corollary \ref{gtlc} of geodesic transversal lemma, we have $gt(BF(k))\geq gt(G[V_1])+gt(G[V_2])\geq 2^{(k-1)+1}+2^{(k-1)+1}=2^{k+1}$, hence $gt(BF(r))\geq 2^{r+1}$, $r\geq 3$ and this completes the proof. 
\end{proof}

\begin{rem}
    For $r=1$, $BF(1)$ is the cycle $C_4$ and the set $S=\{(0,0), (1,1)\}$ is a gt-set. Hence, $gt(BF(1))=2$.
\end{rem}
\begin{rem}
    For $r=2$, clearly $S=\{(00,1),(01,1),(10,1),(11,1)\}$ is a gt-set and hence $gt(BF(2))=4$.
\end{rem}
\section{Hypercube}
A $r$-dimensional hypercube network is a popular interconnection network, with attractive properties such as recursive structure, maximum connectivity, vertex and edge symmetry, etc.

\begin{figure}[htbp]
\centering
\begin{tikzpicture}[scale=0.75,
    x=0.01cm,
    y=-0.01cm,
    vertex/.append style={minimum size=2mm}
]

\node[vertex] (v101a) at (100,45) {};
\node[vertex] (v111)  at (752,45) {};

\node[vertex] (v100)  at (276,225) {};
\node[vertex] (v101b) at (574,225) {};

\node[vertex] (v000)  at (276,523) {};
\node[vertex] (v011)  at (574,523) {};

\node[vertex] (v010)  at (100,701) {};
\node[vertex] (v001)  at (752,701) {};

\draw[rededge] (v101a)--(v111);
\draw[rededge] (v101a)--(v010);
\draw[rededge] (v111)--(v001);
\draw[rededge] (v010)--(v001);

\draw[green!70!black, line width=1pt]
    (v100)--(v101b);
\draw[green!70!black, line width=1pt]
    (v101b)--(v011);
\draw[green!70!black, line width=1pt]
    (v011)--(v000);
\draw[green!70!black, line width=1pt]
    (v000)--(v100);

\draw[edge] (v101a)--(v100);
\draw[edge] (v111)--(v101b);
\draw[edge] (v010)--(v000);
\draw[edge] (v001)--(v011);

\node[anchor=east] at (88,45) {$101$};
\node[anchor=west] at (766,45) {$111$};

\node[anchor=west] at (286,255) {$100$};
\node[anchor=east] at (564,255) {$101$};

\node[anchor=west] at (286,493) {$000$};
\node[anchor=east] at (564,493) {$011$};

\node[anchor=east] at (88,701) {$010$};
\node[anchor=west] at (766,701) {$001$};

\end{tikzpicture}
    \caption{$Q_2'$ and $Q_2''$ in $Q_3$}
    \label{fig6:hypercybe}
\end{figure}
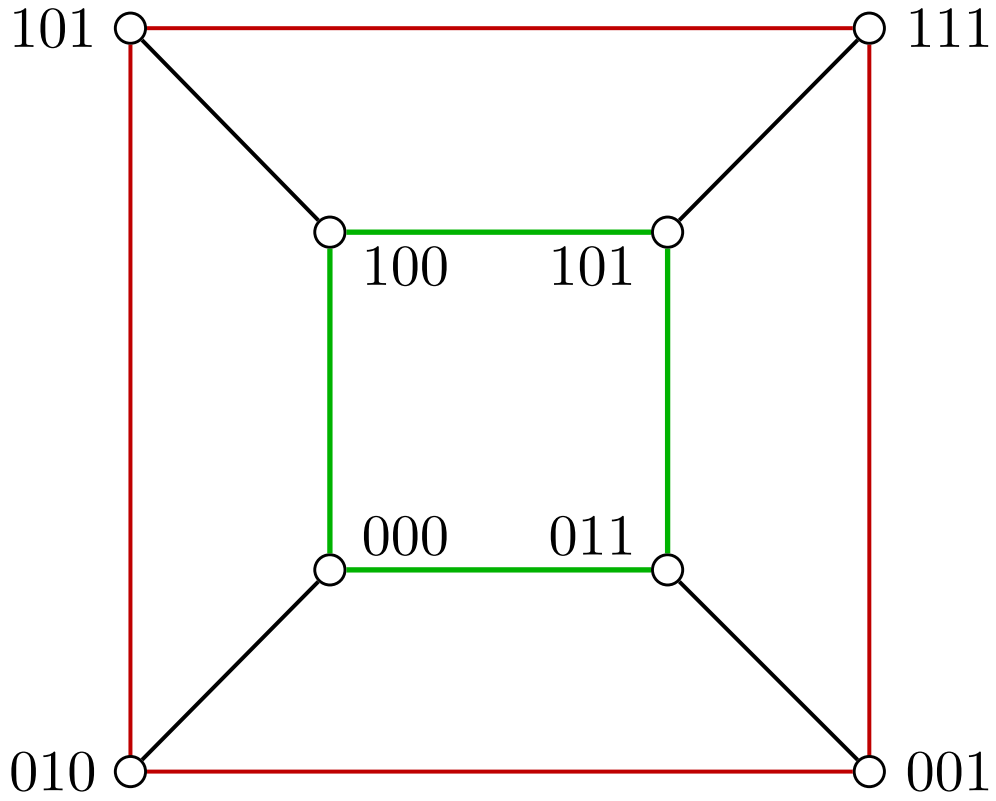

For $r\geq 1$, let $Q_r$ represent the $r$-dimensional hypercube. The vertex set of $Q_r$ is formed by the collection of all $r$-dimensional binary representations. Two vertices $x,y\in G$ are connected by an edge if and only if their binary strings vary precisely by one bit\cite{grid.Qn}. The diameter of $Q_r$ is $r$ and we observe that $Q_r$ contains two hypercubes of dimension $r-1$, which we denote as $Q_{r-1}'$ and $Q_{r-1}''$, refer Figure~\ref{fig6:hypercybe}.

\begin{thm}
    If $r> 1$, then $gt(Q_r)=2^{r-1}$.
\end{thm}
\begin{proof}
   For $r=2$, $Q_2$ is the cycle $C_4$ of diameter 2. Clearly, no single vertex can be a geodesic transversal of $Q_2$ and the set of vertices $\{00,11\}$ is a gt-set of $Q_2$. Hence the theorem is true for $r=2$.

   Let us assume that the result is true for $r=k-1$, that is $gt(Q_{k-1})=2^{k-2}$. We claim that $gt(Q_k)=2^{k-1}$. Observe that every maximal geodesic of $Q_k$ is a diametral path of length $k$ and passes through both the copies $Q_{k-1}'$ and $Q_{k-1}''$. Therefore, the vertex set $V(Q_{k-1}')$(respectively $V(Q_{k-1}'')$) is a geodesic transversal of $Q_k$ and hence $gt(Q_k)\leq 2^{k-1}$.

   Let $P_1$ and $P_2$ be maximal geodesics of $Q_{k-1}'$ and $Q_{k-1}''$ respectively. Both $P_1$ and $P_2$ are diametral paths of length $k-1$. Suppose there exists a path $P$ containing $P_1$ and $P_2$ in $Q_k$, then $P$ cannot be a maximal geodesic since the length of $P$ is greater than $2(k-1)>k$. Therefore, no maximal geodesic of $Q_k$ contains two maximal geodesics either both from $Q_{k-1}'$(respectively $Q_{k-1}''$) or one from each of $Q_{k-1}'$ and $Q_{k-1}''$ respectively. Now let $G[V_1]$ and $G[V_2]$ be $Q_{k-1}'$ and $Q_{k-1}''$ respectively and by our hypothesis $gt(G[V_1])=gt(Q_{k-1}')=2^{k-2}$ and $gt(G[V_2])=gt(Q_{k-1}'')=2^{k-2}$. Applying Corollary \ref{gtlc} of geodesic transversal lemma, $gt(G)=gt(Q_k)\geq gt(G[V_1])+gt(G[V_2])\geq 2^{k-2} + 2^{k-2}=2^{k-1}$ and this completes the proof. 
\end{proof}

\section{Conclusion}
The study of Geodesic Transversal problem is significant because of its strong connectivity with classical concepts like vertex covers, dominating sets and geodesic convexity. Due to its distance constraints, it becomes more relevant in communication networks, transportation systems and social networks. The key result of our paper namely, the Geodesic Transversal Lemma is of great significance in the study of the Geodesic Transversal problem as the gt-number of recursive graph architectures with known induced subgraph decompositions can be determined using this lemma. By applying this lemma, we have solved the geodesic transversal problem of Butterfly networks and Hypercube networks.

\end{document}